\documentclass{amsart}

\usepackage{amsmath}
\usepackage{amssymb}
\usepackage{amsthm}
\usepackage{fullpage}
\usepackage{hyperref}
\usepackage{color}

\font\eighttt=cmtt8

\def\1{{\overline{1}}}
\def\2{{\overline{2}}}
\def\frac#1#2{{#1 \over #2}}
\begin{document}

\title{\bf How To Gamble If You're in a Hurry}
\author[S.B.E.]{Shalosh B. Ekhad}

\author[D.Z.]{Doron Zeilberger}
\thanks{
D.Z. was supported in part by the USA National Science Foundation.
}

\date{\today}

\title{\bf How To Gamble If You're in a Hurry}
\maketitle

{\bf Preamble}

This article is a brief description of the 
Maple package~\href{http://www.math.rutgers.edu/\~zeilberg/tokhniot/HIMURIM}{
{\tt HIMURIM}} downloadable from  

\href{http://www.math.rutgers.edu/\~zeilberg/tokhniot/HIMURIM}{\tt http://www.math.rutgers.edu/\~{}zeilberg/tokhniot/HIMURIM} .

Sample input and output files can be obtained from the webpage:

\href{http://www.math.rutgers.edu/\~zeilberg/mamarim/mamarimhtml/himurim.html}{\tt http://www.math.rutgers.edu/\~{}zeilberg/mamarim/mamarimhtml/himurim.html} .

The Maple package~\href{http://www.math.rutgers.edu/\~zeilberg/tokhniot/HIMURIM}{
{\tt HIMURIM}} is to be considered as the {\bf main output} of this project, and the present article is to be considered as a short user's manual.

We also briefly describe another Maple package downloadable from  

\href{http://www.math.rutgers.edu/\~zeilberg/tokhniot/PURIM}{\tt http://www.math.rutgers.edu/\~{}zeilberg/tokhniot/PURIM} .

{\bf How To Gamble If You Must}

Suppose that you currently have $x$ dollars, and you enter a casino with the hope of getting
out with $N$ dollars, (with, $x$ and $N$, being positive integral values). The probability of winning {\it one} round is $p$ ($0<p<1$).
You can stake any  integral amount of dollars $s(x)$ (that must satisfy $ 1 \leq s(x) \leq min(x,N-x)$), 
until you either exit the casino with the hoped-for $N$ dollars, or you become broke.
Deciding the value of the stake $s(x)$, for each $1 \leq x <N$, constitutes your {\it strategy}. Naturally, the question of whether a strategy is {\it optimal} arises; the three main optimality criteria in gambling theory can be summarized as follows: 
\begin{itemize}
\item[{\bf 1}] Maximizing the probability of 
reaching a specified goal (i.e., amount $N$), with no time limit.
\item[{\bf 2}] Maximizing the probability of reaching a specified goal by a fixed time $T$.
\item[{\bf 3}] Minimizing the expected time to reach a specified goal, subject to a pre-specified level
of {\it risk-aversion}. 
\end{itemize}

In their celebrated masterpiece,  Dubins and Savage~\cite{DS} proved that
the optimal strategy (using the first criterion), if $p \leq \frac{1}{2}$,
is the {\bf bold} one taking $s(x)=min(x,N-x)$, always betting the maximum,
and   if $p \geq 1/2,$ then an optimal strategy is the {\it timid} one, with
$s(x)=1$, always betting the minimum.

A beautiful, lucid, and accessible account of these results can be found in Kyle Siegrist's~\cite{S} on-line article.

Alas, if you play timidly, i.e. according to the classical ``gambler's ruin'' problem (~\cite{F}, p. 348, Eq. (3.4))
your {\it expected time} until exiting is (let $q:=1-p$)

\[
\left\{\begin{array}{ll}
\frac{x}{q-p}-\frac{N}{q-p} \frac{1-(q/p)^x}{1-(q/p)^N} & \textrm{if $p \neq \frac{1}{2};$}\\
x(N-x),& \textrm{if $p=\frac{1}{2},$}
\end{array}\right.
\]

and this may take a very long time. If $p >\frac{1}{2}$, but you're in a hurry, then you may decide 
to take a slightly higher chance of exiting as a loser if that will  enable you to
expect to leave the casino much sooner. It turns out that the bold strategy is way too risky.
For example, if $p=3/5$ and right now you have $100$ dollars and the exit amount is $200$ dollars,
with the bold strategy, sure enough, you are guaranteed to exit the casino after just one round,
but your chance of leaving as a winner is only $3/5$.

As a compromise,  we can employ a {\it deterministic fixed fractional} betting strategy, namely, the {\it Kelly strategy}\footnote{The Kelly strategy is also known as the {\it Kelly system} or the {\it Kelly criterion}; terms first coined by Ed Thorp in~\cite{Thorp1962} and~\cite{Thorp1966}. The theoretical underpinnings of this strategy were provided by Breiman in~\cite{B}.}, with factor $f$ denoting a fixed fraction of our money. This is inspired by~\cite{K}, 
 however, in that paper, the underlying assumptions are: money is infinitely divisible, the game continues indefinitely, the game has even payoff, and the opponent is infinitely wealthy. Under those circumstances, Kelly recommends to take $f=2p-1$ for his agenda. Based on our set of assumptions -- using integral values -- we obtain,
$$
K(f)(x):=min(\lfloor x f \rfloor +1 , N-x) \quad .
$$
For example, the Kelly strategy
with $f=1/10$ (and still $p=3/5,x=100,N=200$) 
enables you to exit as a winner with probability
$\%99.98784517$, but the expected duration is only $44.94509484$ rounds, much sooner than 
the expected duration of $500$ rounds (with a fat tail!) promised by the timid strategy. 

Inspired by Breiman~\cite{B} we can generalize the Kelly-type strategy, and ``morph'' it with the bold strategy, and play boldly once our capital is $\geq Nc$, in other words

\[
B(f,c)(x):= \left\{\begin{array}{ll}
min(\lfloor x f \rfloor + 1, N-x), & \textrm{if $x<cN;$}\\
min(x,N-x), & \textrm{if $x \geq cN.$}
\end{array}\right.
\]

For example, taking $f=1/10, c=4/5$ 
(and still $p=3/5,x=100,N=200$),
your probability of exiting as a winner is
$\%99.98721302$, only slightly less than Kelly with $f=1/10$, but your expected stay at the
casino is about one round less ($43.81842784$). Paradoxically, lowering the $c$ to $2/5$ is not advisable, since
your probability of winning is lower ($\%99.93836900$) {\it and} you should expect to stay longer!
($52.61769977$ rounds). We observed, empirically, that for any $f$, lowering the $c$ from $1$ until a certain place $c_0(f)$ reduces the expected duration-until-winning (with a slightly higher risk of ultimate loss), but setting $c$ below $c_0$ (i.e., playing boldly starting at $cN$) will not only lower your chance of ultimately winning, but would also prolong your agony of staying in the casino (unless you want to {\it maximize} your stay there, in which case you should play timidly).

Our question is: what is the optimal strategy
according to Criterion 2 (i.e. maximizing the probability of reaching a specified goal by a fixed time $T$)?
Borrowing the colorful yet gruesome language of~\cite{CSYZ}, you owe $N$ dollars to a loan shark who would kill you if you don't return the debt in $\leq T$ units time (rounds of gambling). Luckily, you are at a 
{\it superfair} casino, 
(i.e. $p\geq \frac{1}{2}$). If your current capital is 
$i$ dollars (so you need to make $N-i$ additional dollars in $\leq T$ 
rounds to stay alive), if you want to {\it maximize} your chances of 
staying alive, how many dollars should you stake ?

{\bf Finding the Best Strategy If You're in a Hurry}

So suppose that you currently have $i$ dollars, and you need to 
make $N-i$ additional dollars,
so that you can exit the casino with $N$ dollars
in $\leq T$ rounds of gambling, where at each  round you can
stake any amount between $1$ and $min(i,N-i)$. 
You want to {\it maximize} your chance of success.
How much should you stake, and what is the resulting probability, let's call it $f(i,T)$.
The probability of winning a single round is $p$.

Obviously $f(i,T)$ satisfies the {\it dynamical programming} recurrence
$$
f(i,T)=
\max \{ (1-p)f(i-x,T-1)+pf(i+x,T-1) \, : \, 1 \leq x \leq \min(i,N-i) \}, \quad
$$
with the obvious {\it boundary conditions} $f(0,T)=0,f(N,T)=1$ and 
{\it initial conditions}
$f(N,0)=1$, and $f(x,0)=0$ if $x<N$.

The set of $x$'s that attain this max constitutes 
your {\it optimal strategy}. It is most convenient to take the largest $x$
(in case there are ties).

By repeatedly computing $f(j,T)$ and the stake-amount $x$ that
realizes it, where $j$ is the current capital and $T$ is
the steadily decreasing time left, the gambler can always
know how much to stake in order to maximize his chance of 
staying alive, and also
know the actual value of that probability.

{\bf The Maple Package HIMURIM}

\href{http://www.math.rutgers.edu/\~zeilberg/tokhniot/HIMURIM}{\tt HIMURIM} is downloadable, {\it free of charge}, from

\href{http://www.math.rutgers.edu/\~zeilberg/tokhniot/HIMURIM}{\eighttt http://www.math.rutgers.edu/\~{}zeilberg/tokhniot/HIMURIM} .

We will only briefly describe some of the more important procedures, leaving it to the readers
to explore and experiment with the many features on their own, using the on-line help.

{\bf The most Important Procedures of HIMURIM}

The most important procedure is  {\tt BestStake(p,i,N,T)}, that implements 
$f(i,T)$ with the given $p$ and $N$.

If you want so see the {\it full} optimal strategy, a list of length $N-1$ whose $i$-th
entry tells you how much to stake if you have $i$ dollars, use
procedure {\tt BestStrat(p,N,T)}. 
See, for example,

~\href{http://www.math.rutgers.edu/\~zeilberg/tokhniot/oHIMURIMk1}{\tt http://www.math.rutgers.edu/\~{}zeilberg/tokhniot/oHIMURIMk1} \quad .

for the output of  {\tt BestStrat(11/20,1000,30);} .

Procedure
{\tt SimulateBSv(p,i,N,T)} simulates {\it one} game that follows 
the optimal strategy.

Finally, Procedure
{\tt BestStratStory(m0,N0,T0,K)} collects optimal strategies for
various $p$'s (all superfair), $N$'s (exit capitals) and $T$ (deadlines).

{\bf Other Procedures of HIMURIM}

Procedure {\tt ezraLA()} lists the procedures that use {\it Linear Algebra} to find the {\it exact} probabilities, expected duration, and
probability generating functions for the random variables ``duration'' and ``duration conditioned on ultimately winning''
for a casino with exit capital $N$, probability of winning a round $p$ (that may be either numeric or {\it symbolic}), and
{\it all} possible initial incomes.

For example, {\tt PrW(p,S)} inputs a probability {\tt p} (a  number between $0$ and $1$ or left as a symbol $p$) and
a list {\tt S}, of length $N-1$, say, where $S[i]$ tells you how much to stake if you have $i$ dollars.
It outputs the list, let's call it {\tt L}, such that {\tt L[i]} is the probability of ultimately winning
(exiting with $N$ dollars)
if you currently have {\tt i} dollars and always play according to strategy {\tt S}.

It works by solving the system  of $N-1$ equations for the $N-1$ unknowns $L[1], \dots, L[N-1]$
$$
L[i]=(1-p)L[i-S[i]]+pL[i+S[i]] \quad ,  1 \leq i \leq N-1 \quad ,
$$
together with the boundary conditions $L[0]=0,L[N]=1$.

For example, if $N=3$,  and the strategy $S$ being $[1,1]$,
(the only possible strategy when $N=3$), then

{\tt PrW(1/3,[1,1]);}

would yield 

{\tt [1/7, 3/7]} ,

that means that if the probability of winning a round is $\frac{1}{3}$ and you exit the casino  when you either reach $0$ or $3$ dollars, then
your probability of exiting as a winner, if you currently have one dollar is $\frac{1}{7}$, and if you currently have $2$ dollars, 
is $\frac{3}{7}$.

Procedure {\tt ED(p,S)} inputs $p$ and $S$ as above
and outputs the list, let's call it {\tt L}, such that {\tt L[i]} is the expected duration until getting out
(either as a winner or loser)
if you currently have {\tt i} dollars and always play according to strategy {\tt S}.

It works by solving the system  of $N-1$ equations for the $N-1$ unknowns $L[1], \dots, L[N-1]$
$$
L[i]=(1-p)L[i-S[i]]+pL[i+S[i]]+1 \quad ,  \quad 1 \leq i \leq N-1 \quad ,
$$
together with the boundary conditions $L[0]=0,L[N]=0$.

For example, still with $N=3$  and $S=[1,1]$,

{\tt ED(1/3,[1,1]);}

would yield 

{\tt  [12/7, 15/7]} ,

that means that if the probability of winning a single round is $\frac{1}{3}$ and you exit the casino  when you either reach $0$ or $3$ dollars, 
and you follow strategy $[1,1]$,
then the expected remaining duration, if you currently have one dollar, is $\frac{12}{7}$ , and if you currently have $2$ dollars, 
it is $\frac{15}{7}$.

Procedure {\tt EDw(p,S)} inputs $p$ and $S$ as above
and outputs the list, let's call it {\tt L}, such that {\tt L[i]} is the expected duration until getting out,
conditioned on being an ultimate winner!
if you currently have {\tt i} dollars and always play according to strategy {\tt S}.

For example, with the above (trivial) input

{\tt EDw(1/3,[1,1]);}

would yield 

{\tt  [18/7, 11/7]} ,

that means that if the probability of winning a round is $\frac{1}{3}$ and you exit the casino  when you either reach $0$ or $3$ dollars, then
the expected remaining duration until winning (assuming that you do win), if you currently have one dollar is 
$\frac{18}{7}$, and if you currently have $2$ dollars is $\frac{11}{7}$.

Procedure {\tt Dpgf(p,S,t)} inputs $p$ and $S$ as above, as well as a {\it symbol} (variable name) $t$,
and outputs the list, let's call it {\tt L}, such that {\tt L[i]} is the 
probability generating function, in $t$, for the random variable ``remaining duration'' if you currently have $i$ dollars,
i.e. if you take the Maclaurin expansion of {\tt L[i]} and extract the coefficient of $t^j$ you would get
the exact value of the probability that the game would last 
exactly $j$ more rounds.

For example, still with the same $N$ and $S$,

{\tt lprint(Dpgf(1/3,[1,1],t));}

would yield 

{\tt [-t*(6+t)/(-9+2*t**2), -t*(3+4*t)/(-9+2*t**2)]} \quad .

Typing ``{\tt taylor(-t*(6+t)/(-9+2*t**2),t=0,5);}'' would yield

{\tt 2/3*t+1/9*t**2+4/27*t**3+2/81*t**4+O(t**5)}

meaning that if you play the above game with $p=1/3,N=3$ and you currently have one dollar,
you would have probability $2/3$ of exiting after one round, 
probability $1/9$ of exiting after two rounds, 
probability $4/27$ of exiting after three rounds, 
and probability $2/81$ of exiting after four rounds.

Procedure {\tt DpgfW(p,S,t)} is analogous to {\tt Dpgf(p,S,t)} but the duration is conditioned on the fortunate event
of exiting as an ultimate winner.
For example, for the timid strategy, and $N=3$, 

{\tt lprint(DpgfW(1/3,[1,1],t));}

would yield 

{\tt  [-t**2/(-9+2*t**2), -3*t/(-9+2*t**2)]} \quad .

Typing {\tt taylor(-t**2/(-9+2*t**2),t=0,5);} would yield

{\tt  1/9*t**2+2/81*t**4+O(t**6)}

meaning that if you play the above game with $p=1/3,N=3$ and you currently have one dollar,
and you are destined to leave as a winner,
you would have probability of $0$ of exiting after one round, 
probability of $1/9$ of exiting after two rounds, 
probability of $0$ of exiting after three rounds, 
probability of $2/81$ of exiting after four rounds.

{\bf The Best Breiman-Kelly Strategies If You are in A Hurry}

Procedure {\tt BestBKdd(p,N,i,T,h)} tells you the best 
Breiman-Kelly Strategy if
the probability of winning a round is {\tt p}, you
have {\tt i} dollars and you must exit the casino with
{\tt N} dollars in $<=$ {\tt T} rounds, and you're
using resolution {\tt h}.
It also returns the expected duration until exit (either as a winner or loser). 

To get the story for various initial capitals, and various deadlines, try out procedure  \hfill\break
{\tt BestBKddStory(p,N,h,t0,MaxF,M0)}. See the on-line help, and the
sample input and output files in

\href{http://www.math.rutgers.edu/\~zeilberg/mamarim/mamarimhtml/himurim.html}{\eighttt http://www.math.rutgers.edu/\~{}zeilberg/mamarim/mamarimhtml/himurim.html} .

{\bf The Best Kelly Factor With a Given Level of Risk-Aversion} 

Try out procedure {\tt KellyContestx(p,N,x,h,conf)}. For example see:

\href{http://www.math.rutgers.edu/\~zeilberg/tokhniot/oHIMURIM8a}{\tt http://www.math.rutgers.edu/\~{}zeilberg/tokhniot/oHIMURIM8a} ,

\href{http://www.math.rutgers.edu/\~zeilberg/tokhniot/oHIMURIM8b}{\tt http://www.math.rutgers.edu/\~{}zeilberg/tokhniot/oHIMURIM8b} .

{\bf The Maple package} {\tt PURIM}

For an ``umbral'' approach, see our Maple package {\tt PURIM} that tells you much more.
It explores the whole ``tree'' of possibilities. See the package itself and the input and output files

\href{http://www.math.rutgers.edu/\~zeilberg/tokhniot/inPURIM2}{\tt http://www.math.rutgers.edu/\~{}zeilberg/tokhniot/inPURIM2} and 

\href{http://www.math.rutgers.edu/\~zeilberg/tokhniot/oPURIM2}{\tt http://www.math.rutgers.edu/\~{}zeilberg/tokhniot/oPURIM2} for an example.

{\bf Further Work}

There are many possible generalizations and extensions.
See for example the interesting article~\cite{CSYZ}. 

{\bf Conclusions}

The three authors completely agree on the {\it mathematics}, but they have somewhat different
views about the significance of this project. Here they are.

{\bf Doron Zeilberger's Conclusion}

Traditional mathematicians like Dubins and Savage 
use traditional proof-based mathematics, 
and also work in the framework of {\it continuous}
probability theory using the pernicious Kolmogorov, 
{\it measure-theoretic}, paradigm. This approach was
fine when we didn't have computers, but we can do so much more 
with both {\it symbol-crunching}
and {\it number-crunching}, in addition to naive simulation, 
and develop {\it algorithms}
and write {\it software}, that  ultimately is a much more useful 
(and rewarding) activity than ``proving''
yet-another-theorem in an artificial and 
fictional continuous, measure-theoretic, world, that
is furthermore {\it utterly boring}.

{\bf Shalosh B. Ekhad's Conclusion}

These humans, they are so emotional! That's why they never went very far.

{\bf Email addresses}: SBE: {\it c/o zeilberg@math.rutgers.edu} ; DZ: {\it zeilberg@math.rutgers.edu} .


\begin{thebibliography}{111}

\bibitem{B} L. Breiman, {\it Optimal gambling systems for favorable games}, Fourth Berkley Symposium on Math. Stat. and Prob. {\bf 1} (1961), 65-78.

\bibitem{BHS} D.~A.~Berry and D.~C.~Heath and W.~D.~Sudderth, {\it Red-and-black with unknown win probability}, Ann. Statist. 2 (1974), 602–608. 

\bibitem{CSYZ} R.W. Chen, L.A. Shepp, Y.-C. Yao, and C.-H. Zhang, {\it On optimality of bold play
for primitive casinos in the presence of inflation}, J. Appl. Prob. {\bf 42} (2005), 121-137.

\bibitem{DS} Lester E. Dubins and Leonard J. Savage, {\it Inequalities for Stochastic Processes}
(How to Gamble If You Must), Dover, 1976.

\bibitem{F} William Feller, ``{\it An Introduction to Probability
Theory and Its Application}'',  Wiley,volume 1, third edition, 1968.

\bibitem{K} J.L. Kelly, Jr, {\it A new interpretation of information rate}, Bell System Tech. J. {\bf 35} (1956), 917-926.

\bibitem{MS} A.P. Maitra and W.D. Sudderth, {\it ``Discrete Gambling and Stochastic Games''}, Springer, 1996.

\bibitem{S}  Kyle Siegrist, {\it How to Gamble if you must}, JOMA v. 8,
~\href{http://www.maa.org/joma/Volume8/Siegrist/RedBlack.pdf}{\tt http://www.maa.org/joma/Volume8/Siegrist/RedBlack.pdf}

\bibitem{Thorp1962} E.O. Thorp, {\it ``Beat the Dealer''}, 2nd Ed., Vintage, New York 1966.

\bibitem{Thorp1966} E.O. Thorp and W. Walden, {\it ``A winning bet in Nevada baccarat.''}, J. Amer. Statist. Assoc., 61 Part I, (1966) 313-328.

\bigskip
\hrule
\bigskip
\end{thebibliography}
\end{document}